\documentclass[12pt,oneside,reqno]{amsart}
\usepackage{mathrsfs}
\usepackage{graphics}
\usepackage{graphicx,color}
\usepackage{amssymb}
\usepackage{subfigure}
\newcommand{\dif}{\mathrm{d}}

\newcommand{\be}{\begin{eqnarray}}
\newcommand{\ee}{\end{eqnarray}}
\newcommand{\ce}{\begin{eqnarray*}}
\newcommand{\de}{\end{eqnarray*}}
\newtheorem{theorem}{Theorem}[section]
\newtheorem{lemma}[theorem]{Lemma}
\newtheorem{remark}[theorem]{Remark}
\newtheorem{definition}[theorem]{Definition}
\newtheorem{proposition}[theorem]{Proposition}
\newtheorem{Example}[theorem]{Example}
\newtheorem{corollary}[theorem]{Corollary}
\def\e{\varepsilon}
\def\t{\theta}
\def\a{\alpha}

\def\[{{\Big[}}
\def\]{{\Big]}}
\def\<{{\langle}}
\def\>{{\rangle}}
\def\({{\Big(}}
\def\){{\Big)}}

\def\no{\nonumber}
\def\bt{\begin{theorem}}
\def\et{\end{theorem}}
\def\bl{\begin{lemma}}
\def\el{\end{lemma}}
\def\br{\begin{remark}}
\def\er{\end{remark}}
\def\bx{\begin{Example}}
\def\ex{\end{Example}}
\def\bd{\begin{definition}}
\def\ed{\end{definition}}
\def\bp{\begin{proposition}}
\def\ep{\end{proposition}}
\def\bc{\begin{corollary}}
\def\ec{\end{corollary}}

\def\cB{{\mathcal B}}
\def\cC{{\mathcal C}}

\def\cM{{\mathcal M}}
\def\cN{{\mathcal N}}

\def\mE{{\mathbb E}}

\def\mP{{\mathbb P}}

\def\mR{{\mathbb R}}

\def\mX{{\mathbb X}}

\def\geq{\geqslant}
\def\leq{\leqslant}

\def\sF{{\mathscr F}}

\def\sX{{\mathscr X}}

\title{Effective Filtering on a Random Slow Manifold*}

\author{Huijie Qiao
\\{School of Mathematics,
Southeast University\\
Nanjing, Jiangsu 211189,  China\\
\MakeLowercase{hjqiaogean@seu.edu.cn}}\\\vskip2mm
Yanjie Zhang\\{Center for Mathematical Sciences \& School of Mathematics and Statistics,
Huazhong University of Science and Technology\\
Wuhan, 430074,  China\\
\MakeLowercase{d201477019@hust.edu.cn}}\\\vskip2mm
Jinqiao Duan\\{Department of Applied Mathematics,
Illinois Institute of Technology, \\
Chicago, IL 60616,USA
\\ \MakeLowercase{duan@iit.edu}}}



\thanks{*This work was supported by NSF of China (No. 11001051, 11371352).}


\date{}

\begin{document}

\allowdisplaybreaks

\begin{abstract}
This work is about a slow-fast data assimilation system under Gaussian noisy fluctuations. First, we obtain its
low dimensional reduction via an invariant slow manifold. Second, we prove that the low dimensional filter   on the slow manifold approximates the original filter in a suitable metric. Finally,  we illustrate this approximate filter numerically   in an   example.
\end{abstract}

\maketitle \rm

\section{Introduction}
  Stochastic dynamical systems evolving on multiple time scales arise widely in   engineering and science.
For example,   dynamics of chemical reaction networks often take place on notably different
time scales, from the order of nanoseconds to the order of several days.  The approximation by two time
  scales is common in various situations. This is especially true for gene regulatory networks (\cite{Tur}),
as the mRNA synthesis process is  significantly faster than the protein dynamics,  and this leads to a two-time-scale system (\cite{wu}).

Treating stochastic differential equations with two-time scales, Khasminskii and Yin (\cite{RG}) developed a
stochastic averaging principle that enables one to average out the fast-varying variables.
The main idea is as follows: under appropriate conditions, with the slow-varying component "fixed," if the
fast-varying component has a stationary distribution, it can be shown that the process represented by the
slow-changing component converges weakly to a limit averaging system.

For random dynamical systems generated by stochastic differential equations with two-time scales, the
theory of invariant manifolds provides another approach for   qualitative analysis of dynamical behaviors, as   invariant manifolds  are geometric
structures to describe  or reduce stochastic dynamics (\cite{Fu, WR, Sch}). Under suitable conditions,
Fu-Liu-Duan (\cite{Fu}) obtained low dimensional reduction of stochastic evolutionary equations with two-time scales
via random slow invariant manifolds.

Filtering   is a procedure to extract system state information with the help
of observations (see \cite{da}). The state evolution and the observations are usually under
noisy  fluctuations. The general idea is to achieve the best estimate for the true system
state, given only   noisy observations for the system. It provides an algorithm for estimating a
signal or state of a random dynamical system based on noisy measurements. Stochastic filtering  is
    important in many practical applications,  from inertial guidance of aircrafts
and spacecrafts to weather and climate prediction. Filtering problems for systems with
two-time scales have  been   studied, with help of stochastic averaging (see \cite{Imkeller, Park1, Park2, Park3} and references
therein).

The goal of this present  paper is to investigate   filtering for stochastic differential equations with
slow and fast time scales.  First, we  obtain a
low-dimensional reduced system on a random  slow manifold,  as in \cite{Fu}.   Then we show that the  filter     of the low
dimensional reduced system  converges to the original filter in an appropriate sense, and this will be numerically illustrated in an example.  As far as we know, this is the first study of filtering on
random slow manifolds.

It is worthwhile to mention that our assumption conditions and method  are different from those in available literature on nonlinear filtering problems for stochastic differential equations with slow and fast time scales. On one hand, existence of random slow manifolds need some special conditions. On the other hand, on random slow manifolds these original stochastic differential equations have no Markov property. That means that some techniques, such as the Zakai equations in \cite{Park1, Park2, Park3} and backward stochastic differential equations in \cite{Imkeller}, do not work. Therefore, we make use of an exponential martingale technique to deal with these nonlinear filtering problems.


\bigskip

This paper is organized as follows. In Section \ref{pre},  we recall basic concepts in random dynamical systems and
random invariant manifolds. The framework for our method for reduced filtering is presented in Section \ref{fra}. In Section \ref{filter}, we present
the nonlinear filtering problem and prove the approximation theorem of the filtering. And then a specific
example is tested  to illustrate our method in Section \ref{expe}. Finally, we summarize our work in Section \ref{conclusion}.

The following convention will be used throughout the paper: $C$ with or without indices will denote
different positive constants (depending on the indices) whose values may change from one place to
another.

\section{Preliminaries}\label{pre}

In the section, we introduce some notations and basic concepts in random dynamical systems.

\subsection{Notation and terminology}

$\mathscr{B}(\mR^n)$ stands for the Borel $\sigma$-algebra on $\mR^n$ and $\cB(\mR^n)$ is the set of
all real-valued uniformly bounded Borel-measurable functions on $\mR^n$.  Let $\cC(\mR^n)$ denote
the set
of all real-valued continuous functions on $\mR^n$, and $\cC_b^1(\mR^n)$ denote the collection of all
functions
of $\cC(\mR^n)$ which themselves and their first-order derivatives are uniformly bounded. We introduce the following norm for $\phi\in\cC_b^1(\mR^n)$:
\ce
\|\phi\|=\max\limits_{x\in\mR^n}|\phi(x)|+\max\limits_{x\in\mR^n}|\triangledown\phi(x)|,
\de
where $\triangledown$ stands for the gradient operator. Moreover, $\cC_c^\infty(\mR^n)$ is the collection of all members of
$\cC(\mR^n)$ with continuous derivatives of all orders and with compact support.

\subsection{Random dynamical systems (\cite{la})}

\bd
Let $(\Omega,\sF,\mP)$ be a probability space, and $(\theta_t)_{t\in\mR}$ a family of measurable
transformations from $\Omega$ to $\Omega$.
We call $(\Omega,\sF,\mP; (\theta_t)_{t\in\mR})$ a metric dynamical system if for each $t\in\mR$,
$\theta_t$ preserves the probability measure $\mP$, i.e.,
$$
\theta_t^*\mP=\mP,
$$
and for $s,t\in\mR$,
$$
\theta_0=1_\Omega, \quad \theta_{t+s}=\theta_t\circ\theta_s.
$$
\ed
\bd\label{rds} Let $(\mX,\sX)$ be a measurable space. A mapping
\ce
\varphi: \mR\times\Omega\times\mX\mapsto\mX, \quad
(t,\omega,x)\mapsto\varphi_t(\omega,x)
\de
with the following properties is called a measurable random dynamical system (RDS), or in short,  a
cocycle:

(i) Measurability: $\varphi$ is
$\mathscr{B}(\mR)\otimes\sF\otimes\sX/\sX$-measurable,

(ii) Cocycle property: $\varphi(t,\omega)$ is continuous
for $t\in\mR$,  and further satisfies the following conditions
\be
\varphi(0,\omega)&=&id_{\mX},  \label{perfect coc1}\\
\varphi(t+s,\omega)&=&\varphi(t,\theta_s\omega)\circ\varphi(s,\omega), \label{perfect coc2}
\ee
for all $s,t\in\mR$ and $\omega\in\Omega$.
\ed

\subsection{Random invariant manifolds (\cite{Duan2015,Sch})}

Let $\varphi$ be a random dynamical system on the normed space $(\mX, \|\cdot\|_{\mX})$ . Then we
introduce a random invariant manifold with respect to $\varphi$.

A family of nonempty closed sets $\cM=\{\cM(\omega)\}_{\omega\in\Omega}$ is called a {\it random
set} if for every $y\in\mX$, the mapping
$$
\Omega\ni\omega\rightarrow dist(y,\cM(\omega)):=\inf\limits_{x\in \cM(\omega)}\|x-y\|_{\mX}
$$
is measurable. $\cM$ is called (positively) invariant with respect to $\varphi$ if
$$
\varphi(t, \omega, \cM(\omega))\subset \cM(\t_t\omega), \quad t\geq 0, \quad \omega\in\Omega.
$$

In the sequel, we consider random sets defined by a Lipschitz continuous graph. Define a function by $$
\Omega\times\mR^n \ni (\omega,x)\rightarrow H(\omega, x)\in\mR^m
$$
such that for all $\omega\in\Omega$, $H(\omega, x)$ is globally Lipschitzian in $x$  and for any
$x\in\mR^n$, the mapping $\omega\rightarrow H(\omega, x)$ is a random variable.  Then
$$
\cM(\omega) :=\{(x, H(\omega, x))|x\in\mR^n\},
$$
 is a random set (\cite[Lemma 2.1]{Sch}). The invariant random set
$\cM(\omega)$ is called a {\it Lipschitz random invariant manifold}.


\section{Framework}\label{fra}

In the section, we present the framework for our reduction method for stochastic filtering and present some results which will be applied in the following
sections.

Let $\Omega^1=C_0(\mR, \mR^n)$ be the set of continuous functions on $\mR$ with values in $\mR^n$
that are zero at the origin. This set is equipped with the compact-open topology. Let $\sF^1$ be its Borel
$\sigma$-algebra and $\mP^1$ the Wiener measure on $\Omega^1$.  Define
$$
\t_t^1\omega_1(\cdot):=\omega_1(\cdot+t)-\omega_1(t), \quad \omega_1\in\Omega^1, \quad t\in\mR.
$$
Then $(\Omega^1, \sF^1, \mP^1, \t_t^1)$ is a metric dynamical system.  Similarly, we define
  $\Omega^2=C_0(\mR, \mR^m)$ and $\sF^2, \mP^2, \t_t^2$. Thus,
$(\Omega^2, \sF^2, \mP^2, \t_t^2)$ is another metric dynamical system.  Introduce
$$
\Omega:=\Omega^1\times\Omega^2, ~\sF:=\sF^1\times\sF^2, ~\mP:=\mP^1\times\mP^2,~
\t_t:=\t_t^1\times\t_t^2,
$$
and then $(\Omega, \sF, \mP, \t_t)$ is a metric dynamical system that is used in the sequel.

\bigskip

Consider the following stochastic slow-fast system
\be\left\{\begin{array}{l}
\dot{x}^\e=A{x}^\e+f(x^\e,y^\e)+\sigma_1\dot{V}, \\
\dot{y}^\e=\frac{1}{\e}B{y}^\e+\frac{1}{\e}g(x^\e, y^\e)+\frac{\sigma_2}{\sqrt{\e}}\dot{W},
\label{slfasy}
\end{array}
\right.
\ee
where $A$ and $B$ are $n\times n$ and $m\times m$ matrices respectively,  and the interaction
functions $f:\mR^n\times\mR^m\rightarrow\mR^n$ and $g:\mR^n\times\mR^m\rightarrow\mR^m$ are
Borel measurable. Moreover,  $V, W$ are mutually independent standard Brownian motions taking values in $\mR^n$
and $\mR^m$ respectively,   $\sigma_1$ and $\sigma_2$ are nonzero real noise intensities, and $\e$ is a small positive
parameter representing the ratio of the two time scales.  We make the following hypotheses:

$(\bf{H_1})$ There exists a $\gamma_1\geq0$ such that
$$
\|A\|\leq\gamma_1,
$$
where $\|A\|$ stands for the norm of the matrix $A$ such that $|Ax|\leq \|A\||x|$ for every $x\in\mR^n$,
and $A$ has no eigenvalue on the imaginary axis.

$(\bf{H_2})$ There exists a $\gamma_2>0$ such that
$$
(By, y)\leq-\gamma_2|y|^2, \quad y\in\mR^m.
$$

$(\bf{H_3})$ There exists a positive constant $L$ such that for all $(x_i, y_i) \in \mR^n\times\mR^m$
\ce
|f(x_1, y_1)-f(x_2, y_2)|\leq L(|x_1-x_2|+|y_1-y_2|),
\de
and
\ce
|g(x_1, y_1)-g(x_2, y_2)|\leq L(|x_1-x_2|+|y_1-y_2|).
\de

$(\bf{H_4})$
$$
\gamma_2>L.
$$

$(\bf{H_5})$ There exist two positive constants $C_f, C_g$ such that
\ce
\sup\limits_{(x,y)\in \mR^n\times\mR^m}|f(x,y)|=C_f,\\
\sup\limits_{(x,y)\in \mR^n\times\mR^m}|g(x,y)|=C_g.
\de

Under the  assumptions $(\bf{H_3})$ and $(\bf{H_5})$, the system (\ref{slfasy}) has a global unique
solution   $(x^\e(t), y^\e(t))$,  with a given  initial value $(x_0, y_0)$. Define the solution operator
$\varphi^\e_t(x_0, y_0):=(x^\e(t), y^\e(t))$, and then $\varphi^\e$ is a random dynamical system.
Introduce two auxiliary systems
\ce
&&\dif\eta=A\eta\dif t+\sigma_1\dif V, \\
&&\dif\xi^\e=\frac{1}{\e}B\xi^\e\dif t+\frac{\sigma_2}{\sqrt{\e}}\dif W.
\de
So, by \cite[Lemma 3.1]{Sch}, there exist two random variables $\eta, \xi^\e$ such that
$\eta(\theta^1_t\omega_1), \xi^\e(\theta^2_t\omega_2)$ solve two above equations, respectively. Set
\ce
&&\bar{x}^\e:=x^\e-\eta(\theta^1_{\cdot}\omega_1),\\
&&\bar{y}^\e:=y^\e-\xi^\e(\theta^2_{\cdot}\omega_2),
\de
and then $(\bar{x}^\e, \bar{y}^\e)$ satisfy the following system
\ce\left\{\begin{array}{l}
\dot{\bar{x}}^\e=A\bar{x}^\e+f(\bar{x}^\e+\eta(\theta^1_{\cdot}\omega_1),\bar{y}^\e+\xi^\e(\theta^2_{\cdot}\omega_2)),
\\
\dot{\bar{y}}^\e=\frac{1}{\e}B{\bar{y}}^\e+\frac{1}{\e}g(\bar{x}^\e+\eta(\theta^1_{\cdot}\omega_1),\bar{y}^\e+\xi^\e(\theta^2_{\cdot}\omega_2)).\label{slfasy1}
\end{array}
\right.
\de
Moreover, $(\bar{x}^\e, \bar{y}^\e)$ generates a random dynamical system denoted by
$\bar{\varphi}^\e$. The following theorem comes from \cite[Theorem 4.2]{Sch}.

\bt\label{inma} (Random slow manifold) \\
Suppose that $\e>0$ is sufficiently small and $(\bf{H_1})$--$(\bf{H_5})$ are satisfied. Then
$\bar{\varphi}^\e$ has a random invariant manifold
$$
\bar{\cM}^\e(\omega)=\left\{\big(x, H^\e(\omega,x)\big), x\in\mR^n\right\},
$$
where for $\omega\in\Omega$,
$$
\sup\limits_{x_1\neq x_2\in\mR^n}\frac{|H^\e(\omega,x_1)-H^\e(\omega,x_2)|}{|x_1-x_2|}
\leq\frac{2(\gamma_2-\a)}{\gamma_2-\a-L},
$$
and $\a$ is a positive number satisfying $\gamma_2-\a>L$.
\et

Based on the relation between $\varphi^\e$ and $\bar{\varphi}^\e$, it holds that $\varphi^\e$ also has a
random invariant manifold
$$
\cM^\e(\omega)=\left\{\big(x+\eta(\omega_1), H^\e(\omega,x)+\xi^\e(\omega_2)\big),
x\in\mR^n\right\}.
$$
By the same deduction as \cite[Theorem 4.4]{Fu}, we could get a reduction system on $\cM^\e$.

\bt\label{Reduction-Theo} (Reduced system on the random slow manifold) \\
Assume that $\e>0$ is sufficiently small and $(\bf{H_1})$--$(\bf{H_5})$ hold. Then for the system
(\ref{slfasy}), there exists the following reduced low dimensional system on the random slow manifold:
\be\left\{\begin{array}{l}
\dot{\tilde{x}}^\e=A\tilde{x}^\e+f\left(\tilde{x}^\e,\tilde{y}^\e\right)+\sigma_1\dot{V},\\
\tilde{y}^\e=H^\e( \theta_{\cdot}\omega,
\tilde{x}^\e-\eta(\theta_{\cdot}^1\omega_1))+\xi^\e(\theta_{\cdot}^2\omega_2),
\end{array}
\right.
\label{redsys}
\ee
such that for   $t >0$ and almost all $\omega$,
\ce
|z^\e(t, \omega)-\tilde{z}^\e(t, \omega)|&\leq&C_{L, \gamma_2, \a}e^{\frac{-\a
t}{\e}}|z^\e(0)-\tilde{z}^\e(0)|,
\de
where $\tilde{z}^\e(t)=(\tilde{x}^\e(t), \tilde{y}^\e(t))$ is the solution of the low dimensional system (\ref{redsys}) with
the initial value $\tilde{z}^\e(0)=(\tilde{x}_0, \tilde{y}_0)$ and $C_{L, \gamma_2, \a}>0$ is a
constant depending on $L, \gamma_2$ and $\a$.
\end{theorem}

\section{An approximate filter on the slow manifold}\label{filter}

In the section we introduce nonlinear filtering problems for the system (\ref{slfasy}) and the reduced
system (\ref{redsys}) on the random  slow manifold, and then study their relation.

\subsection{Nonlinear filtering problems}

In the subsection we introduce nonlinear filtering problems for the system (\ref{slfasy}) and the reduced system (\ref{redsys}).

For $T>0$, an observation system  is given by
\ce
r^{\e}_t= U_t+\int_0^t h(x^\e_s, y^\e_s) \dif s, \quad t\in[0,T],
\de
where $U$ is a standard Brownian motion independent of $V$ and $W$. For the observation system $r^{\e}$, we make the following additional hypothesis:

\medspace

$(\bf{H_6})$ $h$ is bounded and Lipschitz continuous in $(x,y)$ whose Lipschitz constant is denoted by
$\|h\|_{Lip}$.

\medspace

Under the assumption $(\bf{H_6})$, $r^{\e}$ is well defined.  Denote
$$
(\Lambda_t^\e)^{-1}:=\exp\left\{-\int_0^t h(x_s^\e, y_s^\e)\dif U_s-\frac12\int_0^t|h(x_s^\e, y_s^\e)|^2\dif
s\right\},
$$
and then $(\Lambda_t^\e)^{-1}$ is an exponential martingale under $\mP$. By use of
$(\Lambda^\e_t)^{-1}$, we  can  define a probability measure $\mP^\e$ via
$$
\frac{\dif \mP^\e}{\dif \mP}=(\Lambda^\e_T)^{-1}.
$$
By the Girsanov theorem for Brownian motions, we can obtain that under the probability measure
$\mP^\e$, $r^{\e}$ is a standard Brownian motion.

Rewrite $\Lambda_t^\e$ as
$$
\Lambda_t^\e=\exp\left\{\int_0^t h(x_s^\e, y_s^\e)\dif r^{\e}_s-\frac12\int_0^t|h(x_s^\e, y_s^\e)|^2\dif
s\right\},
$$
and define
 $$
 \rho_t^\e(\phi) :=\mE^\e[\phi(x_t^\e)\Lambda^\e_t|\mathcal{R}_t^\e], \quad \phi\in \cB(\mR^n),
 $$
 where  $\mE^\e$ stands for the expectation under $\mP^\e$, $\mathcal{R}_t^\e \triangleq\sigma(r_s^\e:
 0\leq s \leq t) \vee \cN$ and $\cN$ is the collection of all $\mP$-measure zero sets. Here $ \rho_t^\e$ is
 called nonnormalized filtering of $x_t^\e$ with respect to $\mathcal{R}_t^\e$.  Introduce  the measure-valued
 process
\ce
 \pi_t^\e(\phi) := \mE[\phi(x_t^\e)|\mathcal{R}_t^\e], \quad \phi\in \cB(\mR^n),
 \de
and then by the Kallianpur-Striebel formula it holds that
\ce
\pi^{\e}_t(\phi)=\frac{\rho^{\e}_t(\phi)}{\rho^{\e}_t(1)}.
\de
Moreover, $ \pi_t^\e$ is called normalized filtering of $x_t^\e$ with respect to $\mathcal{R}_t^\e$, or the nonlinear filtering problem for $x_t^\e$ with respect to $\mathcal{R}_t^\e$.

Besides, we rewrite the reduced system (\ref{redsys}) as
\ce
\dot{\tilde{x}}^\e=A\tilde{x}^\e+\tilde{f}^\e(\omega, \tilde{x}^\e)+\sigma_1\dot{V},
\de
where $\tilde{f}^\e(\omega, x):=f(x, H^\e( \theta_{\cdot}\omega,
x-\eta(\theta_{\cdot}^1\omega_1))+\xi^\e(\theta_{\cdot}^2\omega_2))$, and study the nonlinear filtering
problem for $\tilde{x}^\e$. Set
\ce
&&\tilde{h}^\e(\omega, x):=h(x, H^\e( \theta_{\cdot}\omega,
x-\eta(\theta_{\cdot}^1\omega_1))+\xi^\e(\theta_{\cdot}^2\omega_2)),\\
&&\tilde{\Lambda}^\e_t:=\exp\left\{\int_0^t\tilde{h}^\e(\omega,\tilde{x}^\e_s)\dif
r^{\e}_s-\frac12\int_0^t|\tilde{h}^\e(\omega,\tilde{x}^\e_s)|^2\dif s\right\},
\de
and then $\tilde{\Lambda}^\e_t$ is an exponential martingale under $\mP^\e$. Thus, we define the
nonnormalized filtering for $\tilde{x}^\e$ by
$$
\tilde{\rho}_t^\e(\phi) :=\mE^\e[\phi(\tilde{x}^\e_t)\tilde{\Lambda}^\e_t|\mathcal{R}_t^\e].
$$
And set
$$
\tilde{\pi}_t^\e(\phi):=\frac{\tilde{\rho}_t^\e(\phi)}{\tilde{\rho}_t^\e(1)},
$$
and then we will prove that $\tilde{\pi}^\e$ could be understood as the nonlinear filtering problem for
$\tilde{x}^\e$ with respect to $\mathcal{R}_t^\e$.

\subsection{The relation between $\pi^{\e}_t$ and $\tilde{\pi}^\e_t$}

In the subsection we will show that a suitable distance   between $\pi^{\e}_t$ and $\tilde{\pi}^\e_t$ converges to
zero as $\e\rightarrow0$. Let us start with two key lemmas.

\bl\label{es1}
Under $(\bf{H_6})$, there exists a constant $C>0$ such that
$$
\mE\left|\tilde{\rho}^{\e}_t(1)\right|^{-p}<\exp\left\{(2p^2+p+1)CT/2\right\}, \quad t\in[0,T], \quad
p>1.
$$
\el
\begin{proof}
Let us compute $\mE\left|\tilde{\rho}^{\e}_t(1)\right|^{-p}$. By the H\"older inequality, it holds that
$$
\mE\left|\tilde{\rho}^{\e}_t(1)\right|^{-p}=\mE^\e\left|\tilde{\rho}^{\e}_t(1)\right|^{-p}\Lambda^\e_T
\leq(\mE^\e\left|\tilde{\rho}^{\e}_t(1)\right|^{-2p})^{1/2}(\mE^\e(\Lambda^\e_T)^2)^{1/2}.
$$
For $\mE^\e\left|\tilde{\rho}^{\e}_t(1)\right|^{-2p}$, notice that
$\tilde{\rho}^{\e}_t(1)=\mE^\e[\tilde{\Lambda}^\e_t|\mathcal{R}_t^\e]$. And then it follows from the
Jensen inequality that
$$
\mE^\e\left|\tilde{\rho}^{\e}_t(1)\right|^{-2p}=\mE^\e\left|\mE^\e[\tilde{\Lambda}^\e_t|\mathcal{R}_t^\e]\right|^{-2p}\leq\mE^\e\left[\mE^\e[|\tilde{\Lambda}^\e_t|^{-2p}|\mathcal{R}_t^\e]\right]=\mE^\e[|\tilde{\Lambda}^\e_t|^{-2p}].
$$
Thus, the definition of $\tilde{\Lambda}^\e_t$ allows  us to obtain that
\ce
\mE^\e[|\tilde{\Lambda}^\e_t|^{-2p}]&=&\mE^\e\left[\exp\left\{-2p\int_0^t\tilde{h}^\e(\omega,\tilde{x}^\e_s)\dif
r^{\e}_s+\frac{2p}{2}\int_0^t|\tilde{h}^\e(\omega,\tilde{x}^\e_s)|^2\dif s\right\}\right]\\
&=&\mE^\e\Bigg[\exp\left\{-2p\int_0^t\tilde{h}^\e(\omega,\tilde{x}^\e_s)\dif
r^{\e}_s-\frac{4p^2}{2}\int_0^t|\tilde{h}^\e(\omega,\tilde{x}^\e_s)|^2\dif s\right\}\\
&&\bullet\exp\left\{\left(\frac{4p^2}{2}+\frac{2p}{2}\right)\int_0^t|\tilde{h}^\e(\omega,\tilde{x}^\e_s)|^2\dif
s\right\}\Bigg]\\
&\leq&\exp\left\{(2p^2+p)CT\right\}\mE^\e\left[\exp\left\{-2p\int_0^t\tilde{h}^\e(\omega,\tilde{x}^\e_s)\dif
r^{\e}_s-\frac{4p^2}{2}\int_0^t|\tilde{h}^\e(\omega,\tilde{x}^\e_s)|^2\dif s\right\}\right]\\
&=&\exp\left\{(2p^2+p)CT\right\},
\de
where the last step is based on the fact that  $\exp\left\{-2p\int_0^t\tilde{h}^\e(\omega,\tilde{x}^\e_s)\dif
r^{\e}_s-\frac{4p^2}{2}\int_0^t|\tilde{h}^\e(\omega,\tilde{x}^\e_s)|^2\dif s\right\}$ is an exponential
martingale under $\mP^\e$.

Similarly, we know that $\mE^\e(\Lambda^\e_T)^2\leq\exp\left\{CT\right\}$. So, by
simple calculation, it holds that
$\mE\left|\tilde{\rho}^{\e}_t(1)\right|^{-p}\leq\exp\left\{(2p^2+p+1)CT/2\right\}$. The proof is
complete.
\end{proof}

\bl\label{es2}
Assume that $(\bf{H_1})$--$(\bf{H_6})$ are satisfied. Then for $\phi\in \cC^1_b(\mR^n)$,
$$
\mE\left|\rho^{\e}_t(\phi)-\tilde{\rho}^{\e}_t(\phi)\right|^p\leq C\|\phi\|^p
(\mE^\e|z^\e(0)-\tilde{z}^\e(0)|^{4p})^{1/4}(e^{\frac{-2\a tp}{\e}}+\e)^{1/2}, \quad t\in[0,T], \quad
p>1,
$$
where the constant $C>0$ is independent of $\e$.
\el
\begin{proof}
For $\phi\in \cC^1_b(\mR^n)$, it follows from the H\"older inequality that
\ce
\mE\left|\rho^{\e}_t(\phi)-\tilde{\rho}^{\e}_t(\phi)\right|^p&=&\mE^\e\left|\rho^{\e}_t(\phi)-\tilde{\rho}^{\e}_t(\phi)\right|^p
\Lambda^\e_T\leq (\mE^\e\left|\rho^{\e}_t(\phi)-\tilde{\rho}^{\e}_t(\phi)\right|^{2p})^{1/2}
(\mE^\e(\Lambda^\e_T)^2)^{1/2}\\
&\leq&\exp\left\{CT/2\right\}(\mE^\e\left|\rho^{\e}_t(\phi)-\tilde{\rho}^{\e}_t(\phi)\right|^{2p})^{1/2}
.
\de
In the following, we estimate $\mE^\e\left|\rho^{\e}_t(\phi)-\tilde{\rho}^{\e}_t(\phi)\right|^{2p}$. Based
on the definitions of $\rho^{\e}_t(\phi), \tilde{\rho}^{\e}_t(\phi)$ and the Jensen inequality, it holds
that
\be
\mE^\e\left|\rho^{\e}_t(\phi)-\tilde{\rho}^{\e}_t(\phi)\right|^{2p}&=&\mE^\e\left|\mE^\e[\phi(x_t^\e)\Lambda^\e_t|\mathcal{R}_t^\e]-\mE^\e[\phi(\tilde{x}^\e_t)\tilde{\Lambda}^\e_t|\mathcal{R}_t^\e]\right|^{2p}\no\\
&=&\mE^\e\left|\mE^\e[\phi(x_t^\e)\Lambda^\e_t-\phi(\tilde{x}^\e_t)\tilde{\Lambda}^\e_t|\mathcal{R}_t^\e]\right|^{2p}\no\\
&\leq&\mE^\e\left[\mE^\e\left[\left|\phi(x_t^\e)\Lambda^\e_t-\phi(\tilde{x}^\e_t)\tilde{\Lambda}^\e_t\right|^{2p}\bigg|\mathcal{R}_t^\e\right]\right]\no\\
&=&\mE^\e\left[\left|\phi(x_t^\e)\Lambda^\e_t-\phi(\tilde{x}^\e_t)\tilde{\Lambda}^\e_t\right|^{2p}\right]\no\\
&\leq&2^{2p-1}\mE^\e\left[\left|\phi(x_t^\e)\Lambda^\e_t-\phi(\tilde{x}^\e_t)\Lambda^\e_t\right|^{2p}\right]\no\\
&&+2^{2p-1}\mE^\e\left[\left|\phi(\tilde{x}_t^\e)\Lambda^\e_t-\phi(\tilde{x}^\e_t)\tilde{\Lambda}^\e_t\right|^{2p}\right]\no\\
&=:&I_1+I_2.
\label{i1i2}
\ee

First, we deal with $I_1$. By the H\"older inequality, it holds that
\begin{equation}
\label{11}
\begin{array}{rcl}
I_1&\leq&\displaystyle
2^{2p-1}(\mE^\e\left[\left|\phi(x_t^\e)-\phi(\tilde{x}^\e_t)\right|^{4p}\right])^{1/2}
(\mE^\e\left|\Lambda^\e_t\right|^{4p})^{1/2} \\[1ex]
&\leq&2^{2p-1}\|\phi\|^{2p}(\mE^\e\left|x_t^\e-\tilde{x}^\e_t\right|^{4p})^{1/2}\Bigg(\mE^\e\exp\left\{4p\int_0^t
h(x_s^\e, y_s^\e)\dif r^{\e}_s-\frac{(4p)^2}{2}\int_0^t|h(x_s^\e, y_s^\e)|^2\dif s\right\}\\[2ex]
&& \bullet \exp\left\{\frac{(4p)^2}{2}\int_0^t|h(x_s^\e, y_s^\e)|^2\dif s-\frac{4p}{2}\int_0^t|h(x_s^\e, y_s^\e)|^2\dif s\right\}
\Bigg)^{1/2}\\[2ex]
&\leq&2^{2p-1}\|\phi\|^{2p}C^{2p}_{L, \gamma_2, \a}e^{\frac{-2\a
tp}{\e}}(\mE^\e|z^\e(0)-\tilde{z}^\e(0)|^{4p})^{1/2}e^{p(4p-1)CT},
\end{array}
\end{equation}
where the last step is based on Theorem \ref{Reduction-Theo} and the fact that the process\\
$\exp\left\{4p\int_0^t h(x_s^\e, y_s^\e)\dif r^{\e}_s-\frac{(4p)^2}{2}\int_0^t|h(x_s^\e, y_s^\e)|^2\dif
s\right\}$ is an exponential martingale under $\mP^\e$.

Next, for $I_2$, we know that
\ce
I_2\leq2^{2p-1}\|\phi\|^{2p}\mE^\e\left[\left|\Lambda^\e_t-\tilde{\Lambda}^\e_t\right|^{2p}\right].
\de
Note that by the It\^o formula, $\Lambda^\e_t$ and $\tilde{\Lambda}^\e_t$ satisfy the following equations,
respectively,
\ce
\Lambda^\e_t=1+\int_0^t\Lambda^\e_s h(x_s^\e, y_s^\e)\dif r_s^\e, \quad
\tilde{\Lambda}^\e_t=1+\int_0^t\tilde{\Lambda}^\e_s \tilde{h}^\e(\omega,\tilde{x}^\e_s)\dif r_s^\e.
\de
Thus, by BDG inequality and the H\"older inequality it holds that
\ce
\mE^\e\left[\left|\Lambda^\e_t-\tilde{\Lambda}^\e_t\right|^{2p}\right]&=&\mE^\e\left[\left|\int_0^t\left(\Lambda^\e_s
h(x_s^\e, y_s^\e)-\tilde{\Lambda}^\e_s \tilde{h}^\e(\omega,\tilde{x}^\e_s)\right)\dif
r_s^\e\right|^{2p}\right]\\
&\leq&C\mE^\e\left[\int_0^t\left|\Lambda^\e_s h(x_s^\e, y_s^\e)-\tilde{\Lambda}^\e_s
\tilde{h}^\e(\omega,\tilde{x}^\e_s)\right|^2\dif s\right]^{p}\\
&\leq&CT^{p-1}\int_0^t\mE^\e\left|\Lambda^\e_s h(x_s^\e, y_s^\e)-\tilde{\Lambda}^\e_s
\tilde{h}^\e(\omega,\tilde{x}^\e_s)\right|^{2p}\dif s\\
&\leq&2^{2p-1}CT^{p-1}\int_0^t\mE^\e\left|\Lambda^\e_s h(x_s^\e, y_s^\e)-\Lambda^\e_s
\tilde{h}^\e(\omega,\tilde{x}^\e_s)\right|^{2p}\dif s\\
&&+2^{2p-1}CT^{p-1}\int_0^t\mE^\e\left|\Lambda^\e_s
\tilde{h}^\e(\omega,\tilde{x}^\e_s)-\tilde{\Lambda}^\e_s
\tilde{h}^\e(\omega,\tilde{x}^\e_s)\right|^{2p}\dif s\\
&=:&I_{21}+I_{22}.
\de
For $I_{21}$, by the similar deduction to $I_1$ we have
\ce
I_{21}&\leq&2^{2p-1}CT^{p-1}\int_0^t\|h\|_{Lip}^{2p}C^{2p}_{L, \gamma_2, \a}e^{\frac{-2\a
sp}{\e}}(\mE^\e|z^\e(0)-\tilde{z}^\e(0)|^{4p})^{1/2}e^{p(4p-1)CT}\dif s\\
&=&2^{2p-1}CT^{p-1}\|h\|_{Lip}^{2p}C^{2p}_{L, \gamma_2,
\a}(\mE^\e|z^\e(0)-\tilde{z}^\e(0)|^{4p})^{1/2}e^{p(4p-1)CT}\frac{\e}{2\a p}[1-e^{\frac{-2\a
tp}{\e}}].
\de
And for $I_{22}$, it follows from the bounded property of $h$ that
\ce
I_{22}\leq2^{2p-1}CT^{p-1}C^{2p}\int_0^t\mE^\e\left|\Lambda^\e_s-\tilde{\Lambda}^\e_s
\right|^{2p}\dif s.
\de
So,
$$
\mE^\e\left[\left|\Lambda^\e_t-\tilde{\Lambda}^\e_t\right|^{2p}\right]\leq
C\e(\mE^\e|z^\e(0)-\tilde{z}^\e(0)|^{4p})^{1/2}+C\int_0^t\mE^\e\left|\Lambda^\e_s-\tilde{\Lambda}^\e_s
\right|^{2p}\dif s.
$$
The Gronwall inequality leads  us to obtain that
$$
\mE^\e\left[\left|\Lambda^\e_t-\tilde{\Lambda}^\e_t\right|^{2p}\right]\leq
C\e(\mE^\e|z^\e(0)-\tilde{z}^\e(0)|^{4p})^{1/2}.
$$
Furthermore,
\be
I_2\leq2^{2p-1}\|\phi\|^{2p}C\e(\mE^\e|z^\e(0)-\tilde{z}^\e(0)|^{4p})^{1/2}.
\label{i2}
\ee

Finally, combining (\ref{i1i2}) with (\ref{11})and(\ref{i2}), we have that
$$
\mE^\e\left|\rho^{\e}_t(\phi)-\tilde{\rho}^{\e}_t(\phi)\right|^{2p}\leq\|\phi\|^{2p}Ce^{\frac{-2\a
tp}{\e}}(\mE^\e|z^\e(0)-\tilde{z}^\e(0)|^{4p})^{1/2}+\|\phi\|^{2p}C\e(\mE^\e|z^\e(0)-\tilde{z}^\e(0)|^{4p})^{1/2},
$$
and then
$$
\mE\left|\rho^{\e}_t(\phi)-\tilde{\rho}^{\e}_t(\phi)\right|^p\leq C\exp\left\{CT/2\right\}\|\phi\|^p
(\mE^\e|z^\e(0)-\tilde{z}^\e(0)|^{4p})^{1/4}(e^{\frac{-2\a tp}{\e}}+\e)^{1/2}.
$$
This proves the lemma.
\end{proof}

Now, we are ready to state and prove the main result in the paper. First, we give out two concepts used in the proof of Theorem \ref{filcon}.

\bd\label{sepapo}
The set $M\subset\cC_b^1(\mR^n)$ strongly separates points in $\mR^n$
when the convergence $\lim\limits_{n\rightarrow\infty}\phi(x_n)=\phi(x), \forall \phi\in M$, for some $x, x_n\in\mR^n$,
implies  that  $\lim\limits_{n\rightarrow\infty} x_n =x$.
\ed

\bd\label{condet}
The set $N\subset\cC_b^1(\mR^n)$ is convergence determining for the topology of weak convergence of probability measures, if
$\mu_n$ and $\mu$ are probability measures on $\mathscr{B}(\mR^n)$, such that
$\lim\limits_{n\rightarrow\infty} \int_{\mR^n}\phi\,\dif\mu_n = \int_{\mR^n}\phi\,\dif\mu$ for any $\phi\in N$, then $\mu_n$ converges weakly to $\mu$.
\ed

\bt\label{filcon} (Approximation by the reduced filter on slow manifold)\\
Assume the hypotheses  $(\bf{H_1})$--$(\bf{H_6})$ hold.  Then for $p>1$,  $\e$ sufficiently small,  and $t\in[0,T]$,  there exists a positive constant $C$ such that for
$\phi\in \cC_b^1(\mR^n)$
$$
\mE|\pi^{\e}_t(\phi)- \tilde{\pi}_t^\e(\phi)|^p\leq C\|\phi\|^p
(\mE|z^\e(0)-\tilde{z}^\e(0)|^{16p})^{1/16}(e^{\frac{-4\a tp}{\e}}+\e)^{1/4}.
$$
Thus, for the  distance $d(\cdot, \cdot)$  in the space of probability measures that induces the
weak convergence, the following approximation holds:
\ce
\mE [d(\pi^{\e}_t, \tilde{\pi}_t^\e)] \leq
C(\mE|z^\e(0)-\tilde{z}^\e(0)|^{16p})^{\frac{1}{16p}}(e^{\frac{-4\a tp}{\e}}+\e)^{\frac{1}{4p}}.
\de
This means the filter for the low dimensional system on the random slow manifold approximates the original filter in this distance $d(\cdot, \cdot)$.
\et
\begin{proof}
For $\phi\in \cC^1_b(\mR^n)$, it follows from Lemma \ref{es1} and \ref{es2} that
\ce
\mE|\pi^{\e}_t(\phi)-
\tilde{\pi}_t^\e(\phi)|^{p}&=&\mE\left|\frac{\rho^{\e}_t(\phi)-\tilde{\rho}^{\e}_t(\phi)}{\tilde{\rho}^{\e}_t(1)}-\pi^{\e}_t(\phi)\frac{\rho^{\e}_t(1)-\tilde{\rho}^{\e}_t(1)}{\tilde{\rho}^{\e}_t(1)}\right|^{p}\\
&\leq&2^{p-1}\mE\left|\frac{\rho^{\e}_t(\phi)-\tilde{\rho}^{\e}_t(\phi)}{\tilde{\rho}^{\e}_t(1)}\right|^{p}+2^{p-1}\mE\left|\pi^{\e}_t(\phi)\frac{\rho^{\e}_t(1)-\tilde{\rho}^{\e}_t(1)}{\tilde{\rho}^{\e}_t(1)}\right|^{p}\\
&\leq&2^{p-1}\left(\mE\left|\rho^{\e}_t(\phi)-\tilde{\rho}^{\e}_t(\phi)\right|^{2p}\right)^{1/2}\left(\mE\left|\tilde{\rho}^{\e}_t(1)\right|^{-2p}\right)^{1/2}\\
&&+2^{p-1}\|\phi\|^{p}\left(\mE\left|\rho^{\e}_t(1)-\tilde{\rho}^{\e}_t(1)\right|^{2p}\right)^{1/2}\left(\mE\left|\tilde{\rho}^{\e}_t(1)\right|^{-2p}\right)^{1/2}\\
&\leq&C\|\phi\|^p (\mE^\e|z^\e(0)-\tilde{z}^\e(0)|^{8p})^{1/8}(e^{\frac{-4\a tp}{\e}}+\e)^{1/4}.
\de
To complete the proof, we only consider $\mE^\e|z^\e(0)-\tilde{z}^\e(0)|^{8p}$. By the H\"older
inequality, it holds that
\ce
\mE^\e|z^\e(0)-\tilde{z}^\e(0)|^{8p}=\mE|z^\e(0)-\tilde{z}^\e(0)|^{8p}(\Lambda_t^\e)^{-8p}\leq(\mE|z^\e(0)-\tilde{z}^\e(0)|^{16p})^{1/2}(\mE(\Lambda_T^\e)^{-16p})^{1/2}.
\de
By   simple calculations,   we  obtain that
\ce
\mE(\Lambda_T^\e)^{-16p}&=&\mE\left(\exp\left\{-16p\int_0^T h(x_s^\e, y_s^\e)\dif
U_s-\frac{16p}{2}\int_0^t|h(x_s^\e, y_s^\e)|^2\dif s\right\}\right)\\
&=&\mE\Bigg[\left(\exp\left\{-16p\int_0^T h(x_s^\e, y_s^\e)\dif U_s-\frac{(16p)^2}{2}\int_0^t|h(x_s^\e,
y_s^\e)|^2\dif s\right\}\right)\\
&&\bullet \exp\left\{\frac{(16p)^2}{2}\int_0^t|h(x_s^\e, y_s^\e)|^2\dif s-\frac{16p}{2}\int_0^t|h(x_s^\e,
y_s^\e)|^2\dif s\right\}\Bigg]\\
&\leq&\exp\{C\frac{(16p)^2-16p}{2}\},
\de
where the last step is based on the fact that $\exp\left\{-16p\int_0^T h(x_s^\e, y_s^\e)\dif
U_s-\frac{(16p)^2}{2}\int_0^t|h(x_s^\e, y_s^\e)|^2\dif s\right\}$ is an exponential martingale under
$\mP$.
Thus,
$$
\mE^\e|z^\e(0)-\tilde{z}^\e(0)|^{8p}\leq C(\mE|z^\e(0)-\tilde{z}^\e(0)|^{16p})^{1/2}.
$$

Next, we know that there exists a countable algebra $\{\phi_i, i=1, 2, \cdots\}$ of $\cC_b^1(\mR^n)$ that strongly seperates points in $\mR^n$. Thus, it follows from Theorem 3.4.5 in \cite{ek} that $\{\phi_i, i=1, 2, \cdots\}$ is convergence determining for the topology of weak convergence of probability measures. For two probability measures $\mu, \tau$ on $\mathscr{B}(\mR^n)$,  define
\ce
d(\mu, \tau):=\sum\limits_{i=1}^{\infty}\frac{|\int_{\mR^n}\phi_i\,\dif\mu-\int_{\mR^n}\phi_i\,\dif\tau|}{2^i}.
\de
Then $d$ is a distance in the space of probability measures on $\mathscr{B}(\mR^n)$. Since $\{\phi_i, i=1, 2, \cdots\}$ is convergence determining for the topology of weak convergence of probability measures, $d$ induces the weak convergence. The proof is complete.
\end{proof}

\section{Numerical experiments}\label{expe}

In this   section,  we present an example to illustrate our filtering method on a random slow manifold.

Consider the following   slow-fast  stochastic system
\begin{eqnarray}
\begin{cases}\label{sfs}
\dot{x}^{\varepsilon}=x^{\varepsilon}+\frac{1}{4}\sin(y^{\varepsilon})+0.01\dot{{V}},\\
\dot{y}^{\varepsilon}=-\frac{1}{\varepsilon }y^{\varepsilon}+\frac{1}{4\varepsilon
}\cos(x^{\varepsilon})+\frac{1}{\sqrt{\varepsilon }}\dot{W},
\end{cases}
\end{eqnarray}
where $A=1, B=-1, f(x,y)=\frac{1}{4}\sin y$ and $g(x,y)=\frac{1}{4}\cos x$. It is easy to justify that
$A, B, f, g$ satisfy $(\bf{H_1})$--$(\bf{H_5})$ with $\gamma_1=\gamma_2=1,
L=C_f=C_g=\frac{1}{4}$.  Then the system (\ref{sfs}) has a unique solution $(x^{\varepsilon},
y^{\varepsilon})$, which generates a  random dynamical system $\varphi^{\varepsilon}$.

Introduce the following two auxiliary systems
\ce
&&\dif\eta=\eta \dif t +0.01\dif V, \\
&&\dif \xi^{\varepsilon}=-\frac{1}{\varepsilon}\xi^{\varepsilon}\dif t
+\frac{1}{\sqrt{\varepsilon}}\dif W.
\de
Then two   equations have the following stationary solutions, respectively,
\ce
\left\{
\begin{aligned}
\eta(\omega_1)&=-0.01\int_{0}^{\infty}e^{-s}\dif V_s(\omega_1), \\
\xi^{\varepsilon}(\omega_2)&=\frac{1}{\sqrt{\varepsilon}}\int_{-\infty}^{0}
e^{\frac{s}{\varepsilon}}\dif W_s(\omega_2).
\end{aligned}
 \right.
\de
Define
\be
\left\{
\begin{aligned}
\bar{x}^{\varepsilon}_t&:=x^{\varepsilon}_t-\eta(\theta_t^{1}\omega_1),\\
\bar{y}^{\varepsilon}_t&:=y^{\varepsilon}_t-\xi^{\varepsilon}(\theta_t^{2}\omega_2),
\end{aligned}
\right.
\label{rexeye}
\ee
and then $(\bar{x}^\e, \bar{y}^\e)$ solve the following equation
\be\left\{\begin{array}{l}
\dot{\bar{x}}^\e=\bar{x}^\e+\frac{1}{4}\sin\(\bar{y}^\e+\xi^\e(\theta^2_{\cdot}\omega_2)\), ~~\bar{x}_0^\e=x\in\mR,
\\
\dot{\bar{y}}^\e=-\frac{1}{\e}{\bar{y}}^\e+\frac{1}{4\e}\cos\(\bar{x}^\e+\eta(\theta^1_{\cdot}\omega_1)\), ~~\bar{y}_0^\e=y\in\mR.
\end{array}
\right.
\label{traneq}
\ee
Thus, by Theorem 3.1, we can get the following random invariant manifold for
$(\bar{x}^{\varepsilon}, \bar{y}^{\varepsilon})$
\ce
\bar{\mathcal{M}^{\varepsilon}}(\omega)=\big\{\big(x, H^{\varepsilon}(\omega,x)\big), x\in
\mathbb{R}\big\},
\de
where
\ce
H^{\varepsilon}(\omega,x)=\frac{1}{4\varepsilon}\int_{-\infty}^{0}e^{\frac{s}{\varepsilon}}
\cos\big(\bar{x}_s^{\varepsilon}(\omega,x)+\eta(\theta_s^1\omega_1)\big)\dif s.
\de
By (\ref{rexeye}), it holds that $\varphi^{\varepsilon}$ has a random invariant manifold
\ce
\mathcal{M}^{\varepsilon}(\omega)=\big\{\big(x+\eta(\omega_1),
H^{\varepsilon}(\omega,x)+\xi^{\varepsilon}(\omega_2)\big), x\in \mathbb{R}\big\}.
\de
Thus, one can obtain the following reduced one dimensional  system on $\mathcal{M}^{\varepsilon}(\omega)$
\ce
\left\{
\begin{aligned}
\dot{\tilde{x}}^{\varepsilon}&=\tilde{x}^{\varepsilon}+\frac{1}{4}\sin(\tilde{y}^{\varepsilon})
+\sigma_1\dot{V}, \\
\tilde{y}^{\varepsilon}&=H^{\varepsilon}\big(\theta_{\cdot}\omega, \tilde{x}^{\varepsilon}
-\eta(\theta_{\cdot}^1\omega_1)\big)+\xi^{\varepsilon}(\theta_{\cdot}^2\omega_2).
\end{aligned}
\right.
\de

Next, the observation system  is given by
\be
\dif r^{\varepsilon}_t=\arctan(x^{\varepsilon}_t)\dif t+\dif U_{t},
\ee
where $h(x,y)=\arctan(x)$. And $h(x,y)$ satisfies $(\bf{H_6})$.


To facilitate numerical simulation, we make some preparations. First, note that $ H^{\varepsilon}(\omega,x)$
has an approximation $H^{0}(\omega,x)+H^{1}(\omega,x)$ (with error $\mathcal{O}(\varepsilon^2)$)
\ce
\begin{aligned}
H^{0}(\omega,x)&=\int_{-\infty}^{0}e^{s}g\big(x+\eta(\theta_{s\varepsilon}^{1}\omega_1),Y_0(s)
+\xi^{\epsilon}(\theta_{s\varepsilon}^{2}\omega_{2})\big)\dif s \\
&=\frac{1}{4}\int_{-\infty}^{0}e^{s}\cos\big(x+\eta(\theta_{s\varepsilon}^{1}\omega_1)\big)\dif s,
\end{aligned}
\de
and
\ce
H^{1}(\omega,x)&=&\nonumber\int_{-\infty}^0
e^{s}\Big\{g_x \cdot\big[sx+\int_0^s
f\big(x+\eta(\theta_{r\varepsilon}^{1}\omega_1),
Y_0(r)+\xi^{\epsilon}(\theta_{r\varepsilon}^{2}\omega_{2})\big)\,\dif r\big]\\{}&&
+g_y\big(x+\eta(\theta_{r\varepsilon}^{1}\omega_1),
Y_0(s)+\xi^{\epsilon}(\theta_{s\varepsilon}^{2}\omega_{2})\big)Y_1(s)\Big\}\,\dif s \\
&=&\nonumber-\frac{1}{4}\int_{-\infty}^0
e^{s}\Big\{\sin\big(x+\eta(\theta_{s\varepsilon}^{1}\omega_1)\big) \cdot\big[sx+\frac{1}{4}\int_0^s
\sin\big(Y_0(r)+\xi^{\epsilon}(\theta_{r\varepsilon}^{2}\omega_{2})\big)\,dr\big]\Big\}\,\dif s.
\de
Here  $Y_0, Y_{1}$ satisfy the following equations, respectively,
\ce
\left\{
\begin{aligned}
Y_0^{'}(s)&=-Y_0(s)+\frac{1}{4}\cos\big(x+\eta(\theta_{s\varepsilon}^{1}\omega_1)\big),\\
Y_0(0)&=H^{0}(\omega,x)
\end{aligned}
\right.
\de
and
\ce
\left\{
\begin{aligned}
Y_{1}^{'}(s)&=-Y_1(s)-\frac{1}{4}\sin\big(x+\eta(\theta_{s\varepsilon}^{1}\omega_1))
\cdot\big[sx+\frac{1}{4}\int_0^{s}
\sin\big(Y_0(r)+\xi^{\epsilon}(\theta_{r\varepsilon}^{2}\omega_{2})\big)\dif r\big] ,\\
Y_1(0)&=H^{1}(\omega, x).
\end{aligned}
\right.
\de
Second, note  that $\frac{1}{\sqrt{\varepsilon}}W_{t\varepsilon}(\omega)$ is a Brownian
motion.
Hence $\psi_{\varepsilon}:\omega\rightarrow \omega$ is defined implicitly by
$W_t(\psi_{\varepsilon}\omega)=\frac{1}{\sqrt{\varepsilon}}W_{t\varepsilon}(\omega)$.
Thus, after a series of simple calculations, we have
\ce
\left\{
\begin{aligned}
\eta(\theta_{\tau\varepsilon}^1(\omega_1))&=-0.01\sqrt{\varepsilon}\int_{0}^{\infty}e^{-u^{'}}
\dif V_{\frac{u^{'}}{\varepsilon}}(\psi_{\varepsilon}\omega_1),\\
\xi^{\varepsilon}\big(\theta_{\tau\varepsilon}^{2}(\omega_2)\big)&=\int_{-\infty}^{0}e^{u}\dif W_u(\psi_{\varepsilon}\omega_2).
\end{aligned}
\right.
\de
Set
\ce
\eta_1(\psi_{\varepsilon}\omega_1):=-\int_{0}^{\infty}e^{-u}\dif V_{u}(\psi_{\varepsilon}\omega_1).
\de
Then $\eta(\theta_{\tau\varepsilon}^{1}\omega_1)$ is  identically   distributed as
$0.01\eta_1(\psi_{\varepsilon}\omega_1)$.

Now  we apply a particle filtering method (\cite{BC, Park1}) to simulate a nonlinear
filter, which approximate the stochastic process $\pi_t$ with discrete random measures
of the form
\begin{equation}
\sum \limits_{i=1}a_i(t)\delta_{{x_t}^{\varepsilon,i}},
\end{equation}
in other words, with empirical distributions associated with sets of randomly located
particles of stochastic mass $a_1(t)$, ~$a_2(t)$,~ $\cdots$, which have stochastic positions ${x_t}^{\varepsilon,1}$,
~~${x_t}^{\varepsilon,2}$,~~$\cdots$. The particle filtering algorithm can be carried out in  the following steps:\\
Step 1: Initialization \\
For $j=1,2, \cdots, n$ \\
\indent{}Sample ${x}^{\varepsilon}(0),{y}^{\varepsilon}(0)$ from $\pi_{0}$. \\
\indent{}$a_j(0):=1$.\\
end for \\
$ \pi_{0}:=\frac{1}{n}\sum_{j=1}^{n}\delta_{({x}^{\varepsilon}(0),{y}^{\varepsilon}(0))}$ \\
Step 2:  Iteration\\
for l:=0 to m-1 \\
\indent{}for j= to n \\
\indent{}\indent{}Using Euler method to generate the Gaussian random vector $ {x}^{\varepsilon}(t+\frac{\triangle
t}{m})$ and $ {y}^{\varepsilon}(t+\frac{\triangle t}{m})$. \\
\indent{}\indent{}$b_j(t+\frac{\triangle t}{m}):=\arctan({x}^{\varepsilon}(t))(r^{\varepsilon}_{t+\frac{\triangle
t}{m}}-r^{\varepsilon}_{t})-\frac{\triangle t}{2m}||\arctan({x}^{\varepsilon}(t))||^2$ \\
\indent{}\indent{}$a_j(t+\frac{\triangle t}{m}):=a_j(t)exp(b_j(t+\frac{\triangle t}{m}))$ \\
\indent{}end for\\
\indent{}$t:=t+\frac{\triangle t}{m}$ \\
\indent{}$\Sigma(t):=\sum_{j=1}^{n}a_{j}(t)$ \\
\indent{}$\pi_{t}^{n}:=\frac{1}{\Sigma(t)}\sum_{j=1}^{n}\delta_{\big({x}^{\varepsilon}(t),{y}^{\varepsilon}(t)\big)}$.
\\
end for \\
Step 3:  Deterministic resampling  \\
Use the Kitagawa's deterministic resampling algorithm,  as described in \cite{TH}.

\bigskip

For the particle filtering algorithm, we  take:  $ \phi(x)=\frac{10x}{1+x^2}, n=200, m=400, \triangle t=0.02, T=8$.
We will compute
the original filter
$$\pi^{\e}_t(\phi)= \mE[\phi(x_t^\e)|\mathcal{R}_t^\e],$$
the reduced filter
$$ \tilde{\pi}_t^\e(\phi)=\frac{\tilde{\rho}_t^\e(\phi)}{\tilde{\rho}_t^\e(1)}, $$
and the mean-square    error
$$
\mE|\pi^{\e}_t(\phi)- \tilde{\pi}_t^\e(\phi)|^2
$$
and plot in the following figures.

As seen in Figure 1 and Figure 2,  it is clear that if the
initial values of the original slow component and the reduced system are the same, the larger $\e$ is, the
larger the fluctuation of the filtering error is. From Figure 3 to Figure 6  it is found that if the difference for
the initial values of the original slow component and the reduced system becomes larger, the fluctuation of
the filtering error is larger.

\begin{figure}[ht]
\begin{minipage}{0.46\linewidth}
\centerline{\includegraphics[width=1\textwidth]{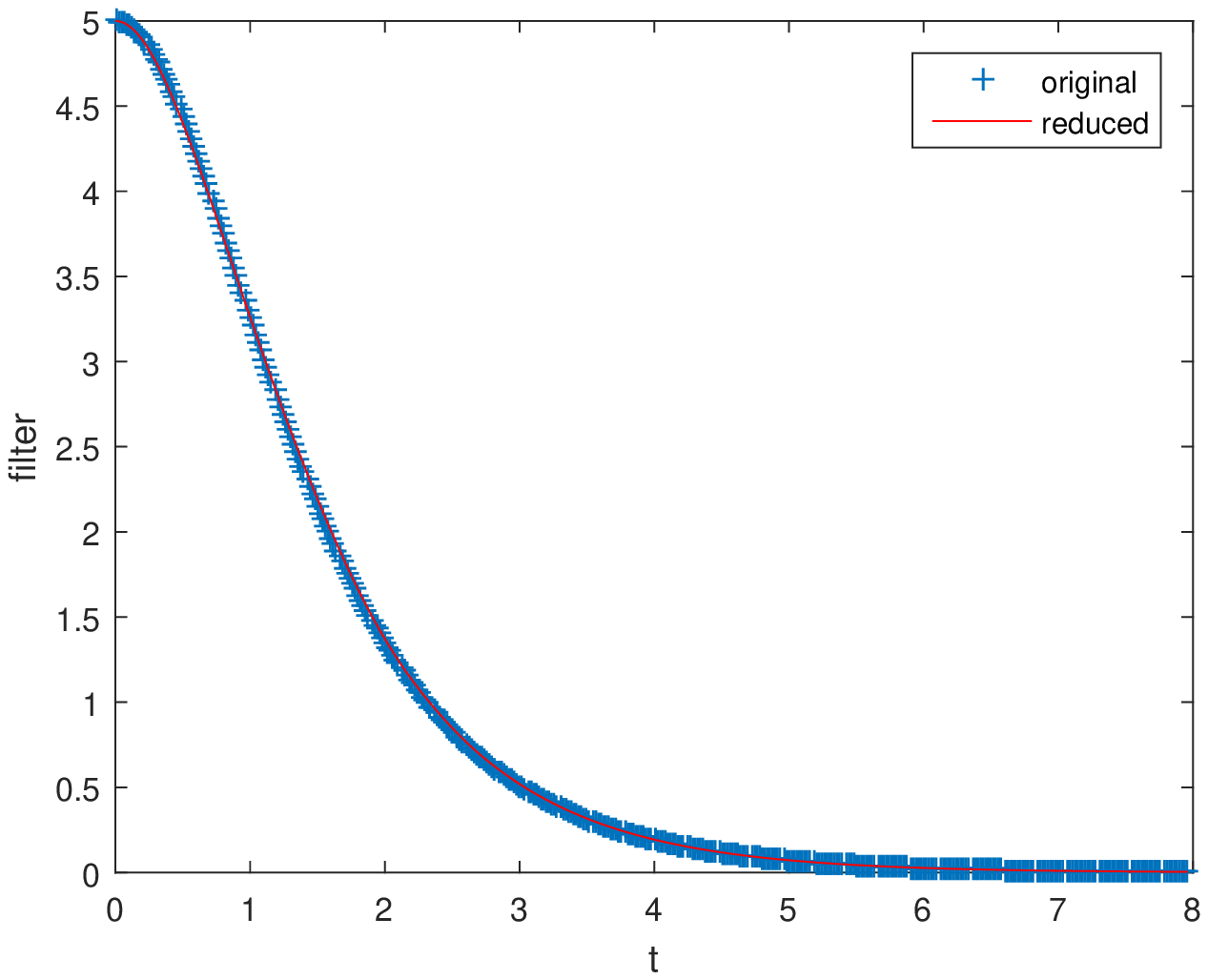}}
\centerline{(a)}
\end{minipage}
\qquad
\begin{minipage}{0.46\linewidth}
\centerline{\includegraphics[width=1\textwidth]{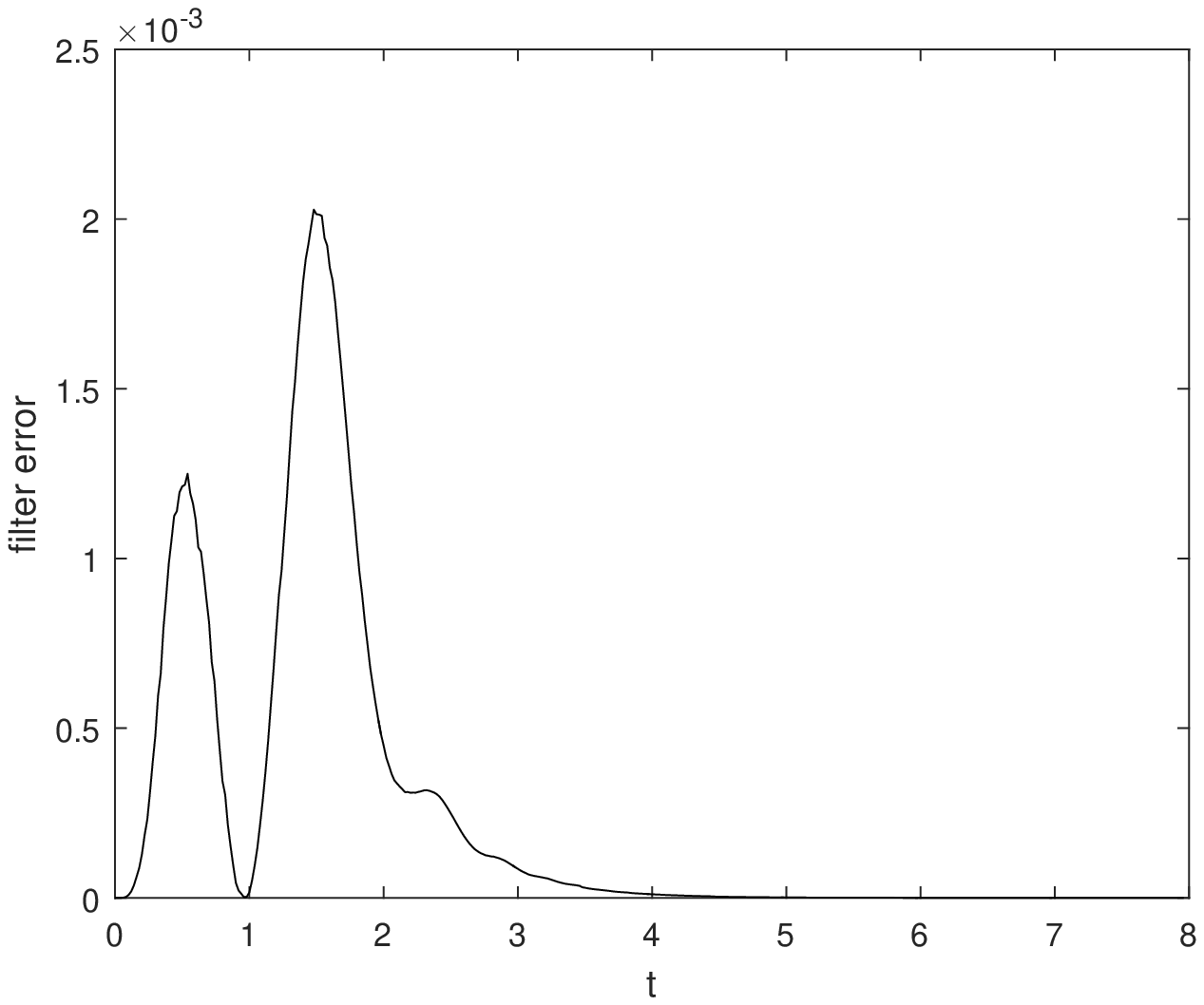}}
\centerline{(b)}
\end{minipage}
\caption{(a) The  original filter $\pi^{\e}_t$   (`+' curves) versus the reduced filter
$ \tilde{\pi}_t^\e$   (red curves):    Initial value $x^{\e}(0)=1$,~~$y^{\e}(0)=1$,~~ $\tilde{x}^{\e}(0)=1$,~$\varepsilon =0.01$;
(b) The mean-square    error  $\mE|\pi^{\e}_t(\phi)- \tilde{\pi}_t^\e(\phi)|^2$.}
\label{Fig1}
\end{figure}

\begin{figure}[ht]
\begin{minipage}{0.46\linewidth}
\centerline{\includegraphics[width=1\textwidth]{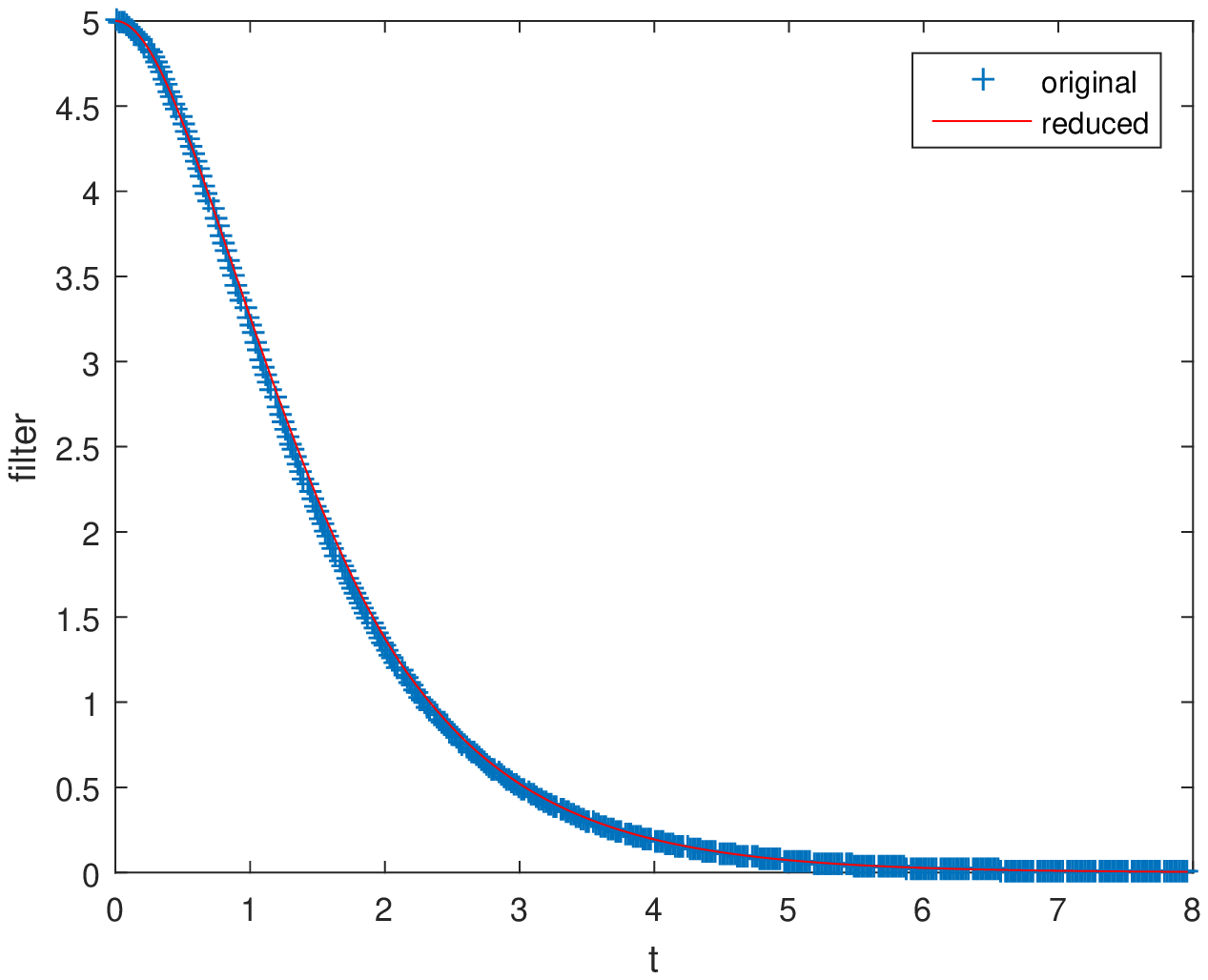}}
\centerline{(c)}
\end{minipage}
\qquad
\begin{minipage}{0.46\linewidth}
\centerline{\includegraphics[width=1\textwidth]{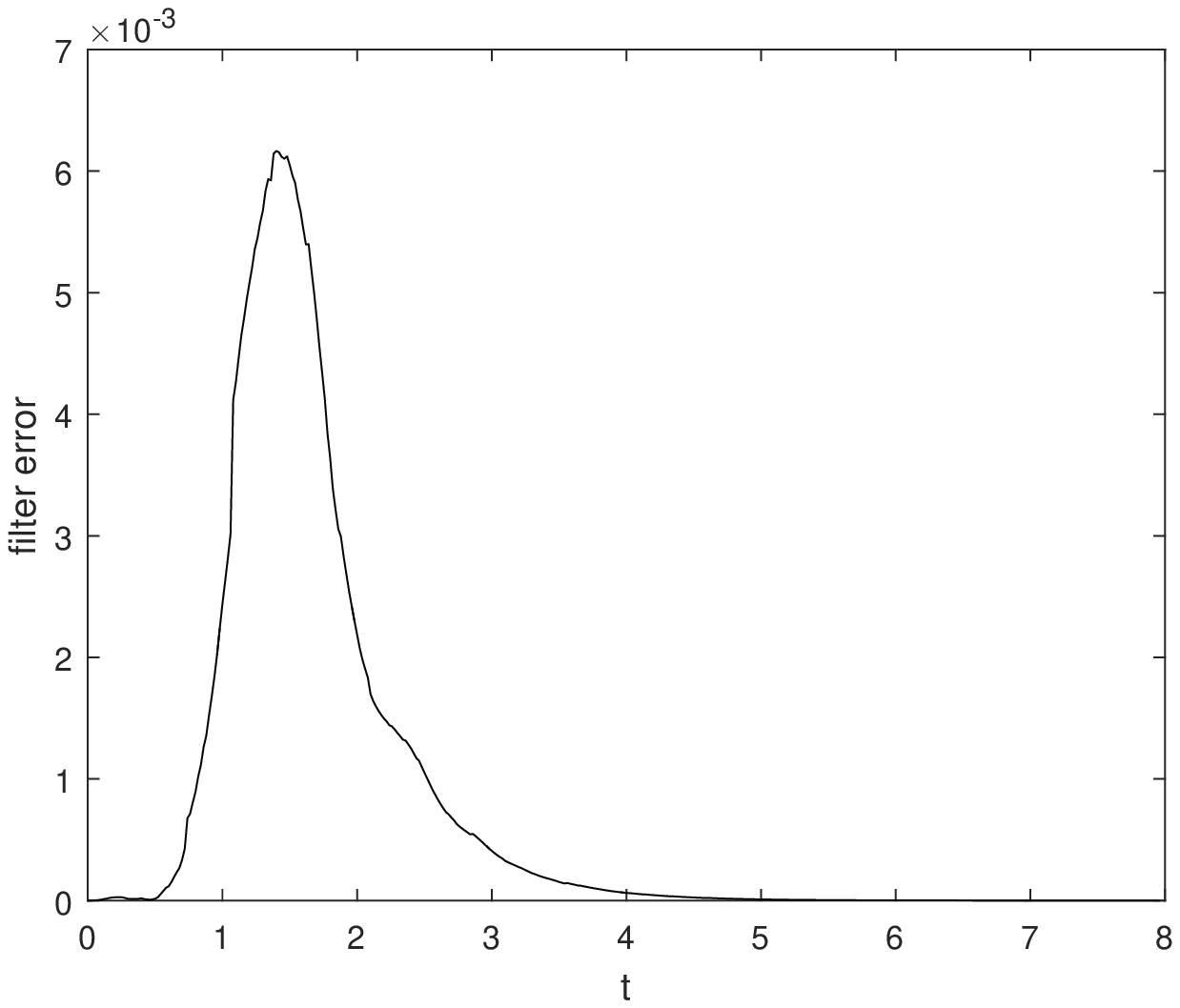}}
\centerline{(d)}
\end{minipage}
\caption{(c) The  original filter $\pi^{\e}_t$   (`+' curves) versus the reduced filter
$ \tilde{\pi}_t^\e$   (red curves): Initial value $x^{\e}(0)=1$,~~$y^{\e}(0)=1$,~~ $\tilde{x}^{\e}(0)=1$,~$\varepsilon =0.1$;
(d) The mean-square error   $\mE|\pi^{\e}_t(\phi)- \tilde{\pi}_t^\e(\phi)|^2$.}
\end{figure}

\begin{figure}[ht]
\begin{minipage}{0.46\linewidth}
\centerline{\includegraphics[width=1\textwidth]{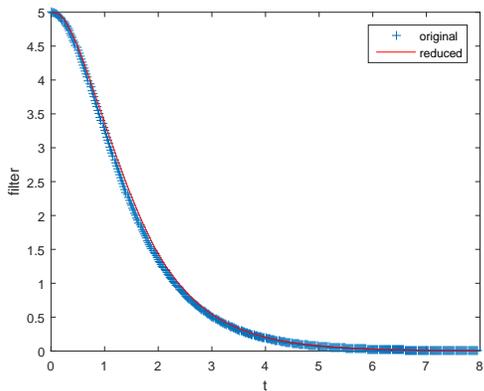}}
\centerline{(e)}
\end{minipage}
\qquad
\begin{minipage}{0.46\linewidth}
\centerline{\includegraphics[width=1\textwidth]{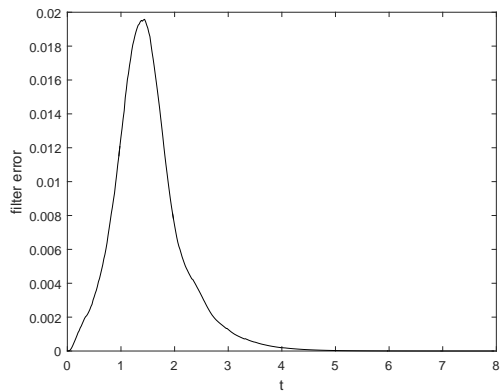}}
\centerline{(f)}
\end{minipage}
\caption{(e) The  original filter $\pi^{\e}_t$   (`+' curves) versus the reduced filter
$ \tilde{\pi}_t^\e$   (red curves):    Initial value $x^{\e}(0)=1$,~~$y^{\e}(0)=1$,~~ $\tilde{x}^{\e}(0)=0.95$,~$\varepsilon =0.01$;
(f) The mean-square error $\mE|\pi^{\e}_t(\phi)- \tilde{\pi}_t^\e(\phi)|^2$.}
\end{figure}

\begin{figure}[ht]
\begin{minipage}{0.46\linewidth}
\centerline{\includegraphics[width=1\textwidth]{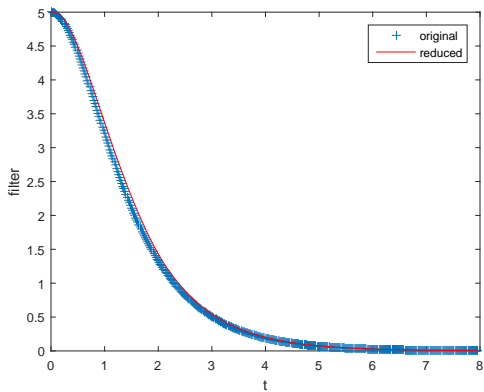}}
\centerline{(g)}
\end{minipage}
\qquad
\begin{minipage}{0.46\linewidth}
\centerline{\includegraphics[width=1\textwidth]{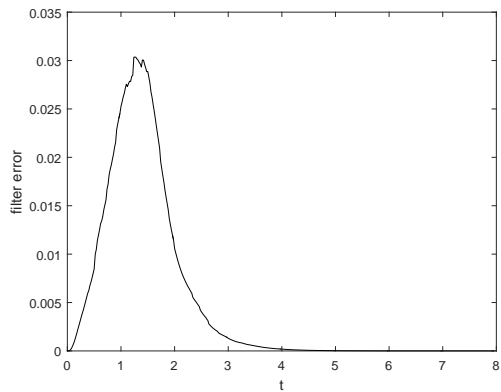}}
\centerline{(h)}
\end{minipage}
\caption{(g) The  original filter $\pi^{\e}_t$   (`+' curves) versus the reduced filter
$ \tilde{\pi}_t^\e$   (red curves): Initial value $x^{\e}(0)=1$,~~$y^{\e}(0)=1$,~~ $\tilde{x}^{\e}(0)=0.95$,~$\varepsilon =0.1$;
(h) The mean-square error  $\mE|\pi^{\e}_t(\phi)- \tilde{\pi}_t^\e(\phi)|^2$.}
\end{figure}

\begin{figure}[ht]
\begin{minipage}{0.46\linewidth}
\centerline{\includegraphics[width=1\textwidth]{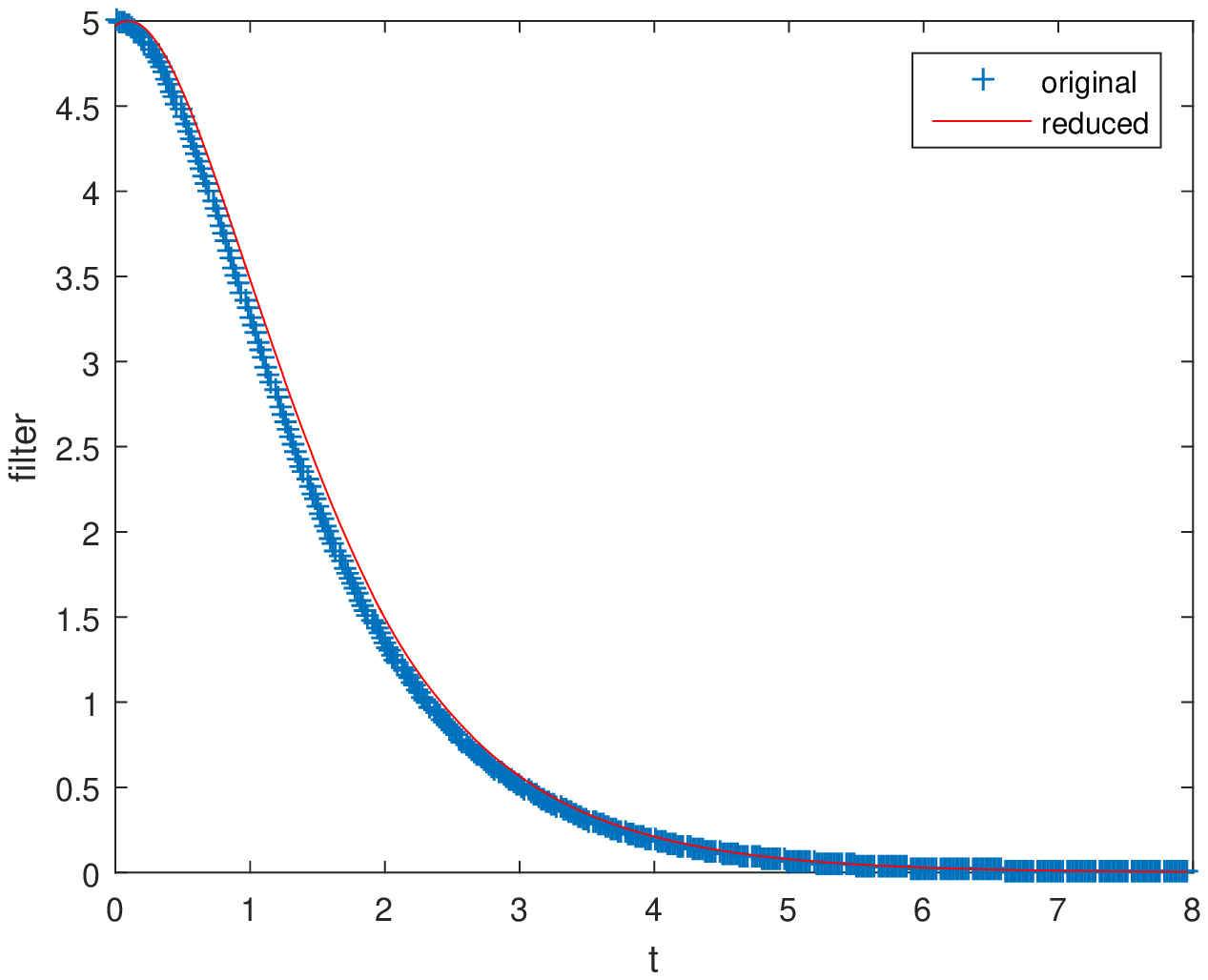}}
\centerline{(i)}
\end{minipage}
\qquad
\begin{minipage}{0.46\linewidth}
\centerline{\includegraphics[width=1\textwidth]{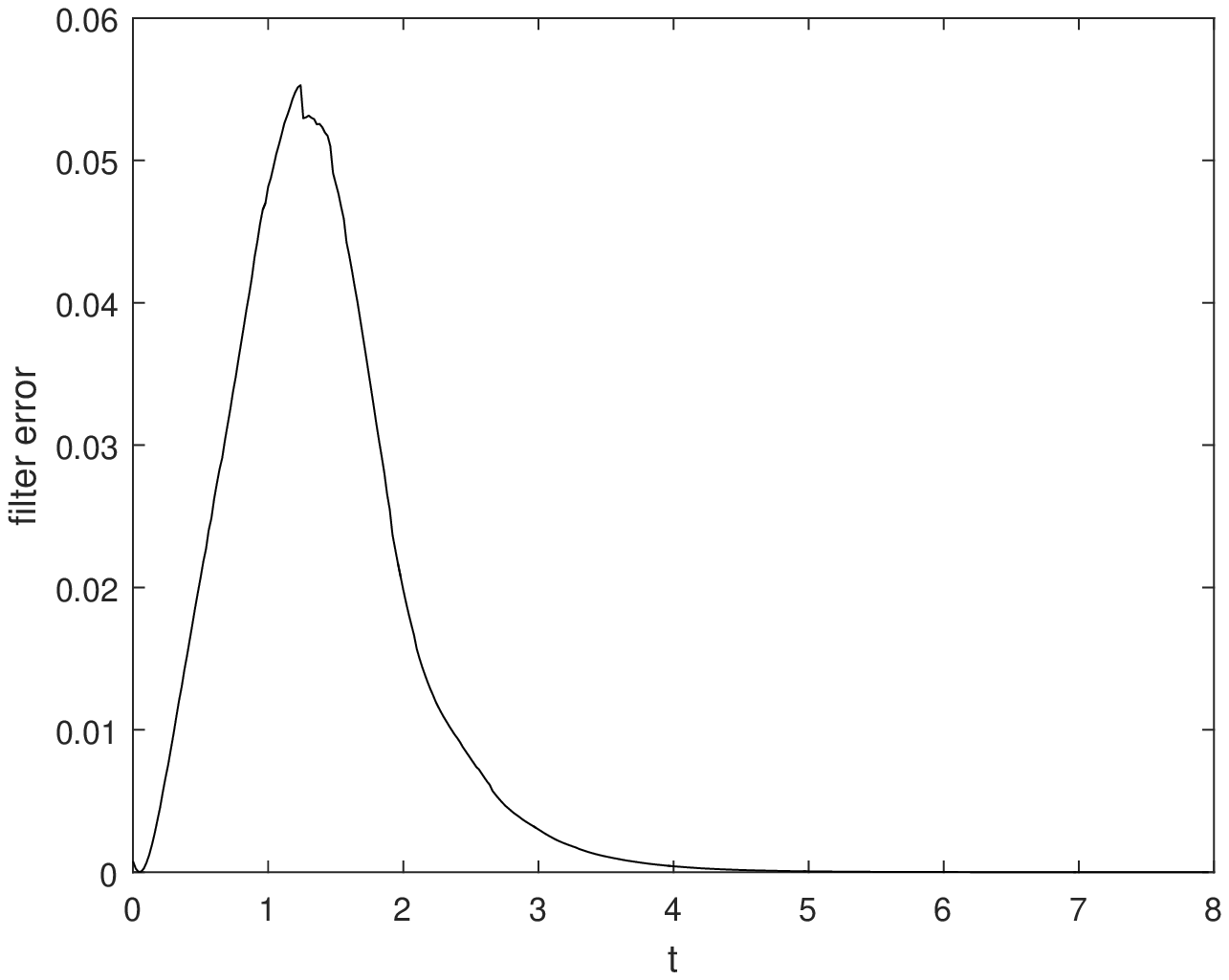}}
\centerline{(j)}
\end{minipage}
\caption{(i) The  original filter $\pi^{\e}_t$   (`+' curves) versus the reduced filter
$ \tilde{\pi}_t^\e$   (red curves):    Initial value $x^{\e}(0)=1$,~~$y^{\e}(0)=1$,~~ $\tilde{x}^{\e}(0)=0.9$,~$\varepsilon =0.01$;
(j) The mean-square error $\mE|\pi^{\e}_t(\phi)- \tilde{\pi}_t^\e(\phi)|^2$.}
\end{figure}
\begin{figure}[ht]
\begin{minipage}{0.46\linewidth}
\centerline{\includegraphics[width=1\textwidth]{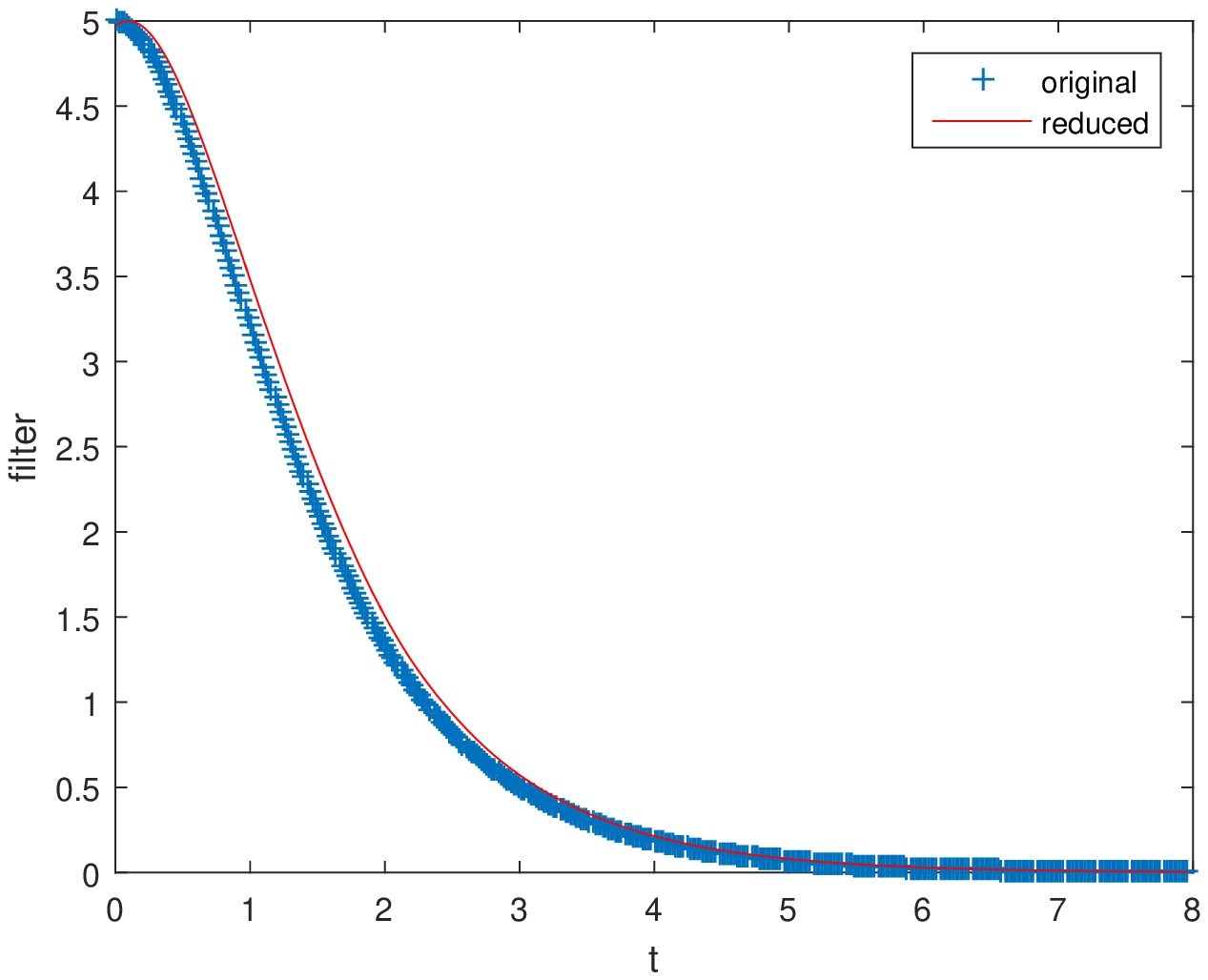}}
\centerline{(k)}
\end{minipage}
\qquad
\begin{minipage}{0.46\linewidth}
\centerline{\includegraphics[width=1\textwidth]{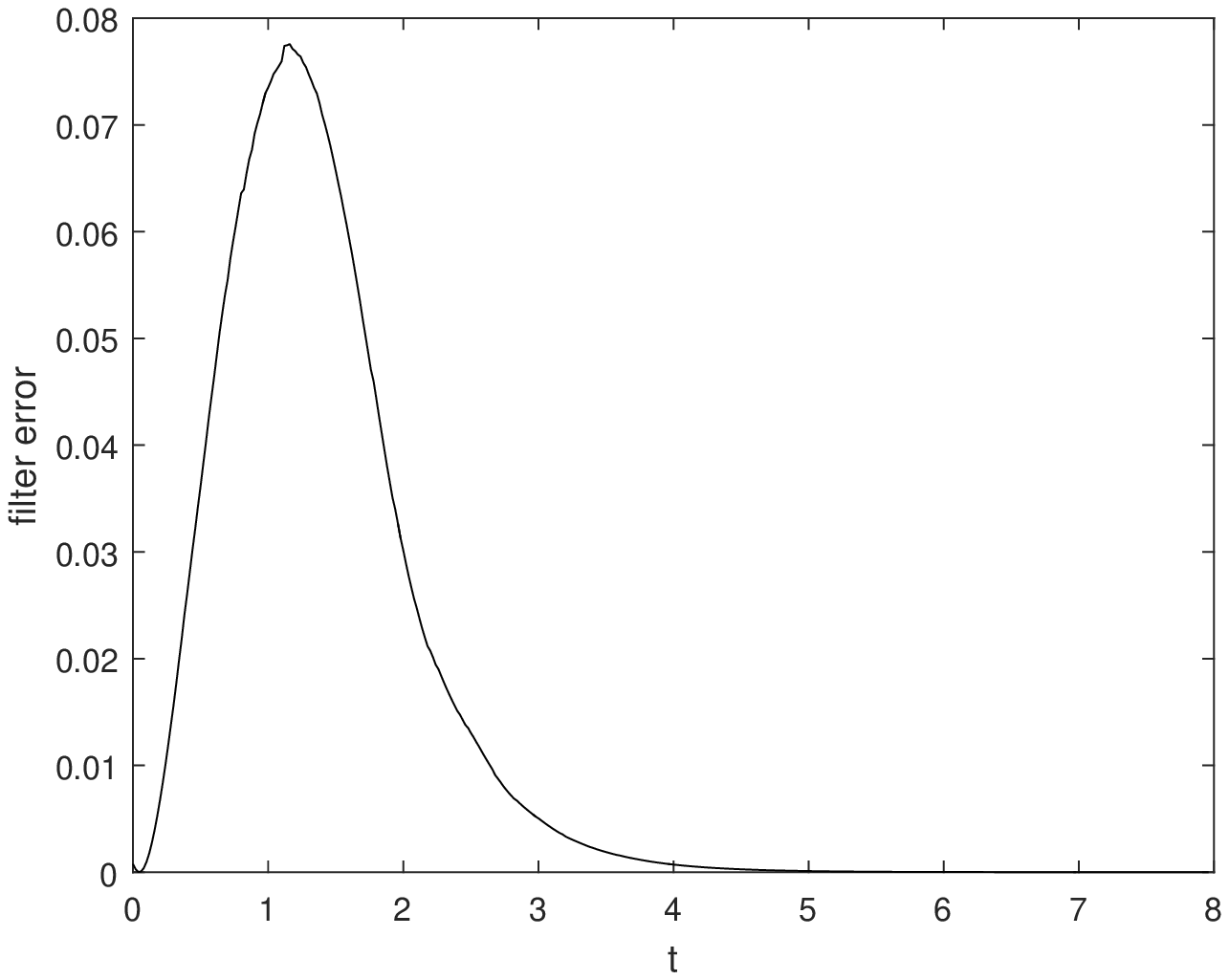}}
\centerline{(l)}
\end{minipage}
\caption{(k)  The  original filter $\pi^{\e}_t$   (`+' curves) versus the reduced filter
$ \tilde{\pi}_t^\e$   (red curves):    Initial value $x^{\e}(0)=1$,~~$y^{\e}(0)=1$,~~ $\tilde{x}^{\e}(0)=0.9$,~$\varepsilon =0.1$;
(l) The mean-square error $\mE|\pi^{\e}_t(\phi)- \tilde{\pi}_t^\e(\phi)|^2$.}
\end{figure}

\section{Conclusion}
\label{conclusion}

In this paper, we first obtain the low dimensional reduction of a slow-fast data assimilation system via an invariant slow manifold.   Then we show that the low dimensional filter on the slow manifold  approximates the original filter. Moreover,  by an example we illustrate this approximate filter numerically.

\bigskip

{\bf Acknowledgements.} We would like to thank Dr. Xinyong Zhang (Tsinghua University, China) and Dr. Jian Ren (Zhengzhou University of Light Industry, China) for helpful
discussions and comments.


\newpage

\begin{thebibliography}{999}
\bibitem{la}
L. Arnold, \emph{Random Dynamical Systems}, Springer, Berlin, 1998.

\bibitem{BC}
A. Bain and D. Crisan, \emph{Fundamentals of Stochastic Filtering}, Springer, Berlin, 2009.

\bibitem{Duan2015} J. Duan,  {\em An Introduction to Stochastic Dynamics}. Cambridge University Press, New York, 2015.

\bibitem{ek} S. N. Ethier and T. G. Kurtz, {\em Markov processes: Characterization and convergence}. John Wiley \& Sons, 1986.

\bibitem{Fu}  H. Fu, X. Liu and J. Duan,
Slow manifolds for multi-time-scale stochastic evolutionary systems,
\emph{Comm. Math. Sci.},
11(2013), 141-162.

\bibitem{TH}
T. Higuchi, On the resampling scheme in the filtering procedure of the Kitagawa Monte Carlo Filter,
Japan, 1995.

\bibitem{Imkeller} P. Imkeller, N. S. Namachchivaya, N. Perkowski and H. C. Yeong,
 Dimensional reduction in nonlinear filtering: a homogenization approach,
\emph{The Annals of Applied Probability},
23(2013), 2290-2326.

\bibitem{Park1}  J. H. Park, N. S. Namachchivaya and H. C. Yeong,
Particle filters in a multiscale environment: Homogenized hybrid particle filter,
\emph{J. Appl. Mech.},
78(2011), 1-10.

\bibitem{Park2}  J. H. Park, R. B. Sowers and N. S. Namachchivaya,
Dimensional reductionin nonlinear filtering,
\emph{ Nonlinearity},
23(2010), 305-324.

\bibitem{Park3}
J. H. Park, B. Rozovskii and R. B. Sowers,
Efficient nonlinear filtering of a singularly perturbed stochastic hybrid system,
\emph{LMS Journal of Computation and Mathematics},
14(2011), 254-270.

\bibitem{Tur}
M. Turcotte,   J. Garcia-Ojalvo,   and G. M. S\"uel,
A Genetic Timer through Noise-Induced Stabilization of an Unstable State,
\emph{Proceedings of the National Academy of Sciences of the United States of America},
41(2008), 15732-7.

\bibitem{RG}
R. Z. Khasminskii and G.Yin,
On transition densities of singularly perturbed diffusions with fast and slow components,
\emph{SIAM J. APPL. MATH.},
56(1996), 1794-1819.


\bibitem{wu}
 F. Wu, T. Tian, J. B. Rawling and George Yin,
 Approximate method for stochastic chemical kinetics with two-time scales by chemical Langevin
 equations,
\emph{Journal of Chemical Physics},
144(2016), 174112.

\bibitem{WR}
W. Wang and A. J. Roberts,
Slow manifold and averaging for slow-fast stochastic differential system,
\emph{J. Math. Anal. Appl},
 398(2013), 822-839.

\bibitem{KP}
Y. Kabanov and S. Pergamenshchikov,
Two-scale stochastic systems: asymptotic analysis and control, Springer-Verlag, Berlin, 2003.


\bibitem{da}
G.  Evensen,
Data Assimilation:The Ensemble Kalman Filter, Springer-Verlag, Berlin, 2009.


\bibitem{RBL}
B. L. Rozovskii,
\emph{Stochastic Evolution System: Linear Theory and Application to nonlinear Filtering},
Springer, New York, 1990.


\bibitem{ZS}
Z.  Schuss,
\emph{Nonlinear Filtering and Optimal Phase Tracking},
Springer, New York, 2012.


\bibitem{par}
E. Pardoux,
Stochastics-An International Journal of Probability and Stochastic Processes,
\emph{Stochastics-An International Journal of Probability and Stochastic Processes},
 3(1979), 127-167.

\bibitem{Sch}
 B. Schmalfu{\ss} and  R. Schneider,
 Invariant manifolds for random dynamical systems with slow and fast variables,
\emph{J. Dyna. Diff. Equa.},
20(2008), 133-164.


\bibitem{andrew}
G. A. Pavliotis and A. M. Stuart,
Multiscale methods: averaging and homogenization,
Springer Science+Business Media, New York, 2008.


\bibitem{ar}
S. Albeverio, B. R\"udiger and J. Wu,
Invariant measures and symmetry property of L\'evy type operators,
\emph{Potential Analysis},
13(2000), 147-168.



\bibitem{abw}
S. Albeverio, Z. Brze\'zniak and J. Wu,
Existence of global solutions and invariant measures for stochastic
differential equations driven by Poisson type noise with non-Lipschitz coefficients,
\emph{ Journal of Mathematical Analysis and Applications},
371(2010), 309-322.






\end{thebibliography}
\end{document}